\documentclass{article}
\usepackage[a4paper, total={6in, 10in}]{geometry}
\usepackage[T1]{fontenc}
\usepackage[utf8]{inputenc}
\usepackage{amsmath,amsthm,amssymb, hyperref}
\usepackage{caption,subcaption,graphicx}
\usepackage{enumerate}
\usepackage{authblk}
\usepackage[sorting=ynt,backend=biber]{biblatex}
\newtheorem{conjecture}{Conjecture}[section]

\title{\textbf{List-Coloring and Chromatic-Choosability - A Dynamic Survey}}
\author{\textbf{Nandana K Vasudevan$^{1}$,  K Somasundaram$^{1}$ and N Narayanan$^{2}$}}
\affil{$^{1}$Department of Mathematics\\ Amrita School of Physical Sciences - Coimbatore\\ Amrita Vishwa Vidyapeetham, India.\\
$^{2}$Department of Mathematics \\
Indian Institute of Technology-Madras, Chennai, India. }
\date{\textit{kv\_nandana@cb.students.amrita.edu, s\_sundaram@cb.amrita.edu, naru@iitm.ac.in}}

\begin{document}

\maketitle

\begin{abstract}
List-coloring, introduced independently by Vizing and by Erd\H{o}s, Rubin, and Taylor in the 1970s, generalizes ordinary vertex coloring by assigning to each vertex its own set of admissible colors.  A graph is chromatic-choosable if its list chromatic number equals its chromatic number. The previous survey on list-coloring  by D R Woodall (2001),  emphasized defective choosability, the list-coloring conjectures, and different methods used for list-coloring. This survey reviews major developments on list-coloring and chromatic-choosability, with emphasis on graph classes for which equality is known, graph classes exhibiting a nontrivial gap, and the principal methods used to prove such results. The survey covers embedded graphs, perfect graphs, complete bipartite and multipartite graphs, claw-free graphs, line graphs, powers of graphs, graph products, and selected variants of list-coloring.
\end{abstract}
\vspace{0.5cm}


\section{Introduction}
\hfill \par The concept of graph list-coloring was introduced independently by Vizing \cite{vizing1976vertex} in 1976  and Erd\H{o}s, Rubin, and Taylor \cite{erdos1979choosability} in 1979. Many foundational results can be found in the seminal paper by Erd{\"os} et al. \cite{erdos1979choosability}. Since list coloring generalizes ordinary vertex coloring, we begin by recalling the basic terminology of vertex coloring.\\

The \textit{proper vertex coloring} or simply \textit{vertex coloring} of a graph is the coloring of its vertices in such a way that no two adjacent vertices receive the same color. We say that a given graph is \textit{$k$-colorable} if the graph admits a proper coloring using $k$ colors. The least $k$ for which a graph $G$ is $k$-colorable is called the \textit{chromatic number} of $G$, denoted by $\chi(G)$.\\

In ordinary vertex coloring, all vertices are colored from a common set of colors.

A list assignment for a graph $G$ is a function $L$ that assigns to each vertex $v$ its palette of colors $L(v)$, called the \textit{list of $v$}. The graph $G$ is said to be  \textit{$L$-colorable} if there exists a proper coloring $c$ of $G$ such that $c(v) \in L(v)$ for every $v \in V(G)$.\\

\par 
A list assignment $L$ is called a $k$-list assignment if $|L(v)|=k$ for every $v\in V(G)$.
 A graph $G$ is \textit{$k$-choosable} or \textit{$k$-list-colorable} if $G$  admits an $L$-coloring for every $k$-list assignment $L$. The least $k$ for which a graph $G$ is $k$-choosable is called the \textit{list chromatic number, choice number}, or \textit{choosability} of $G$, denoted by $\chi_L(G), \chi_\ell(G)$ or $ch(G)$. Throughout this paper, we use the notation $\chi_\ell(G)$.\\

\par  Every $k$-choosable graph is $k$-colorable: take the same $k$-element list at every vertex.  However, any $k$-colorable graph need not be $k$-choosable. For example, any graph $G$ is not $(\chi(G)-1)$-choosable as $G$ is not $L$-colorable when $L(v)=\{1,2,\dots,\chi(G)-1\}$ for all $v\in V(G)$. Here, we find the relation connecting the chromatic number and the list chromatic number. For any graph $G$, the chromatic number acts as a tight lower bound for the list chromatic number, that is, $\chi_\ell(G) \geq \chi(G)$. This inequality is called the \textit{choice chromatic inequality} and was named so by Gravier and Maffray \cite{gravier1998graphs} in 1998.\\

\par 
The equality case $\chi_\ell (G) = \chi(G)$ is the central focus of this survey. Graphs satisfying this equality were called  \textit{$\chi$-choosable graphs} by Gravier and Maffray, while Ohba \cite{ohba2002chromatic}  later used the term \textit{chromatic-choosable graphs}. Basic examples include paths, trees, cycles, and complete graphs. Several major conjectures ask whether broader graph classes are chromatic-choosable; these are discussed in later sections. \\

\par Clearly, not all graphs are chromatic-choosable. The complete bipartite graph $K_{2,4}$ is one such example. We know $\chi(K_{2,4})=2$. Assign the lists $\{\{a,b\},\{c,d\}\}$ and $\{\{a,c\},\{a,d\},\{b,c\},\{b,d\}\}$ for the partite set with two vertices and four vertices, respectively. This is a 2-list assignment for which $K_{2,4}$ is not list-colorable. Thus, $K_{2,4}$ is not 2-choosable, and hence $\chi_\ell(K_{2,4})\neq 2$. For any graph $G$, its  \textit{jump} is defined by $jump(G)= \chi_\ell(G)-\chi(G)$. Erd\H{o}s et al. \cite{erdos1979choosability} observed that the jump can be arbitrarily large. It was demonstrated using the example of the complete bipartite graph $K_{m,m}$, which is not $k$-choosable, where $m=\binom{2k-1}{k}$.\\

\par The remaining part of the survey is organized as follows: Section \ref{Upper bounds} reviews general upper bounds for the list chromatic number, including bounds arising from degeneracy, the coloring number, and the Alon–Tarsi method. Sections \ref{Planar graphs} and \ref{Toroidal graphs} discuss embedded graphs, focusing on planar and toroidal graphs. Sections \ref{Perfect graphs}, \ref{Bipartite graphs} and \ref{Multipartite graphs} cover perfect graphs, complete bipartite graphs, and complete multipartite graphs. Section \ref{Claw-free graphs} concerns claw-free graphs, including claw-free perfect graphs and line graphs. Sections \ref{Square of graphs} and \ref{Powers of graphs} survey results on squares and higher powers of graphs, while Section \ref{Graphs operations} discusses list-coloring under standard graph operations. Section \ref{Some regular graphs} briefly reviews selected results on regular graphs, and Section \ref{Some variants of list-coloring} introduces several variants of list-coloring.

\section{Upper bounds for $\chi_\ell(G)$}\label{Upper bounds}
\subsection{In terms of $\Delta(G)$}
\hfill \par The classical greedy bound gives $ \chi(G) \leq \Delta(G)+1$. Brooks’ theorem improves this to $ \chi(G) \leq \Delta(G)$  for every  connected graph that is neither complete  nor an odd cycle.   Erd\H{o}s et al. \cite{erdos1979choosability} proved the corresponding list-coloring analogue: any connected graph $G$ that is neither a complete graph nor an odd cycle,  is $\Delta(G)$-choosable. The greedy coloring argument also extends to list-coloring, yielding $\Delta (G)+1$  as an upper bound for the list chromatic number for all graphs. 

\subsection{The coloring number}
\hfill \par Another relevant upper bound is the coloring number. The \textit{coloring number} of a graph $G$, denoted by col$(G)$, is given by col$(G)=1+\max \{\delta(G'): G' \subseteq G\}$. Equivalently, col$(G)$, is the least integer $d$ for which $G$ admits an ordering $v_1,v_2,...,v_n$ such that each vertex $v_i$ has at most $d-1$ neighbors among $v_1,v_2,...,v_{i-1}$. Thus col$(G)$ is one more than the degeneracy of $G$, and in particular col$(G) \leq \Delta(G)+1$.

\subsection{The Alon-Tarsi number}
\hfill \par The \textit{Alon-Tarsi number} $AT(G)$ gives a powerful upper bound on the list chromatic number. It originates in the polynomial method of Alon and Tarsi \cite{alon1992colorings}  and was later formalized in the terminology of Jensen and Toft  in their seminal book \cite{jensen2011graph}. Two equivalent viewpoints are commonly used: an algebraic formulation via the graph polynomial and a combinatorial formulation via orientations.\\

\par The first  viewpoint is algebraic and uses the graph polynomial associated with the given graph $G$. This is done by giving a linear order for the vertices and then associating with each vertex $z$, a variable $x_z$. The graph polynomial of $G$ is defined by, 
    \begin{center}
        $f_G(x)=\displaystyle \prod_{uv \in E(G),u<v} (x_u-x_v).$
    \end{center}
If $f_G$ has a non-vanishing monomial $x_g=\prod_{v \in V(G)} x_v^{g(v)}$ with $g(v) \leq f(v)-1$ for $v \in V(G)$, where $f:V(G)\rightarrow \mathbb{N}$, we say that $G$ is \textit{Alon-Tarsi f-choosable} or \textit{f-AT}. In particular, $G$ is \textit{$k$-AT} if $f(v)=k$ for all $v \in V(G)$. Then 
    \begin{center}
        \textit{AT(G)}=$\min \{k:G$ is $k$-AT\}.
    \end{center}

\par The second method relies on the idea of orientation and is combinatorial in nature. The difference of a directed graph $D$, $diff(D)$, can be defined in two ways :
    
\begin{itemize}
    \item $diff(D)= |\varepsilon_{even}(D)|-|\varepsilon_{odd}(D)|$, \\
    where $\varepsilon_{odd}(D)$(or $\varepsilon_{even}(D)$) denotes the set of spanning subdigraphs of $D$ maintaining an equal in-degree and out-degree with an odd(or even) number of edges.
    \item $diff(D)= \sum_{D' \in \mathcal{O}} sign(D')$, \\
    where $\mathcal{O}$ denotes the set of all orientations on the underlying undirected graph $G$ such that for each vertex, the out-degrees are the same as in $D$, where $sign(D')=(-1)^n$, where $n$ is the number of edges that differ in orientation in $D$ and $D'$.  
\end{itemize}
Let $D$ be an orientation on the given graph $G$. An orientation with non-zero difference is called an \textit{Alon-Tarsi orientation}.
    \begin{center}
        $AT(G)=\min\{k:G$ has an $AT-$orientation $D$ with $\Delta^+_D=k-1\}.$
    \end{center}
    
Later, Alon \cite{alon1993restricted} observed that for a graph that can be decomposed into edge-disjoint cliques $Q_1, Q_2,\dots, Q_k$, the set of orientations $\mathcal{O}$ can be narrowed down to the set of orientations $\mathcal{O'}$ that are acyclic on every $Q_i$. Then $diff(D)= \sum_{D' \in \mathcal{O'}} sign(D')$.  Any acyclic orientation is  clearly an AT-orientation.\\ 

\par Despite its strength, the Alon–Tarsi method can be difficult to apply directly. On the algebraic side, dense graphs lead to graph polynomials with many factors and complicated coefficient structure. On the orientation side, the number of relevant orientations grows rapidly with the number of edges. Consequently, the method is most effective when additional structure restricts the possible orientations or simplifies the polynomial.\\

\par To find the coloring number, we have a linear ordering of the vertices. For every $u>v$ with $x_u,x_v$ adjacent in $G$, give the orientation from $x_u$ to $x_v$. This gives an acyclic orientation $D$ of $G$ with $\Delta^+_D \leq $col$(G)-1$. This implies that $AT(G)\leq$ col$(G)$. Thus, summarizing all details together, for any graph $G$, the upper bounds mentioned above occur in the following order:
\begin{center}
$\chi_\ell(G) \leq AT(G)\leq$ col$(G) \leq \Delta(G)+1$.
\end{center}

\subsection{A few other bounds}
\hfill \par We examine a few alternative bounds. Noel et al. \cite{noel2015beyond} established an upper bound given by $\max\{\chi(G), \lceil \frac{\chi(G)+n(G)-1}{3} \rceil\}$. Another bound is in terms of the matching number of the complement of the graph. For any graph $G$, $\chi_\ell(G) \leq n(G)-\bar{\mu}$, where $\bar{\mu}$ is the matching number of $\bar{G}$ \cite{gravier1998graphs}. The best bound that can be acquired using this inequality is $\lfloor \frac{n}{2} \rfloor$, which can be mostly improved. Finally, we can find the list version of Finck's theorem \cite{finck1968chromatic} proved by Erd\H{o}s et al. \cite{erdos1979choosability} which states that for any graph $G$, $\chi_\ell(G)+\chi_\ell(\bar{G}) \leq n(G)+1$. From this inequality, we obtain $\chi_\ell(G)\leq n(G)+1-\chi_\ell(\bar{G})$. This is a helpful bound only when we have the list chromatic number of the complement.\\

\par Few other coloring parameters like paintability number and DP-coloring number also act as upper bounds for list chromatic number. The definition of these concepts and parameter is detailed in Section 13 (13.1 and 13.2). 

\section{Planar graphs}\label{Planar graphs}
\hfill \par A $k$-degenerate graph is $(k+1)$-choosable. Planar graphs being 5-degenerate are 6-choosable. We start with the conjecture posed by Erd\H{o}s et al. \cite{erdos1979choosability} which was the first conjecture on the list-coloring of planar graphs. 

\begin{conjecture}[\cite{erdos1979choosability}] Every planar graph is 5-choosable.
\end{conjecture}

\par Thomassen \cite{thomassen1994every} 
 proved this conjecture by strengthening the induction hypothesis in a highly effective way. He showed that a near-triangulation remains list-colorable when two adjacent vertices on the outer face are precolored, the other boundary vertices have lists of size at least three, and the interior vertices have lists of size at least five.  This carefully designed rigid boundary condition is the key mechanism that allows the induction step to carry through smoothly by carefully peeling a vertex off the boundary and updating the lists of its internal neighbors without dropping their sizes below the required thresholds. Later, the result was broadened to $K_5$-minor-free graphs \cite{skrekovski1998choosability}.\\

\par By the four-color theorem, every planar graph is 4-colorable. Hence it is natural to ask whether every planar graph is also 4-choosable. Erd\H{o}s, Rubin, and Taylor conjectured that this is not the case.

\begin{conjecture}[\cite{erdos1979choosability}]
There exists a planar graph which is not 4-choosable. 
\end{conjecture} 

\par This conjecture was proved true by Voigt \cite{voigt1993list} by giving a non-4-choosable planar graph  with 238 vertices. Later, Gutner \cite{gutner1996complexity} improved this result by giving a smaller graph with 75 vertices that is not 4-choosable. Both of these graphs are 4-colorable but not 4-choosable. Mirzakhani \cite{mirzakhani1996small} subsequently constructed an example of a 3-colorable planar graph with 63 vertices that is not 4-choosable. \\

\par We next survey the sufficient conditions under  which planar graphs are 4-choosable. The initial set of results used forbidden cycles. Lam et al. \cite{lam19994} proved the 4-choosability of $C_4$-free planar graphs. Fijav{\v{z}} et al. \cite{fijavvz2002planar} extended the result to $C_5$-free planar graphs and $C_6$-free planar graphs by showing that they are 3-degenerate. They immediately conjectured the 4-list-colorability of planar graphs without $C_7$. Farzad \cite{farzad2009planar} proved this conjecture and extended the result to planar graphs without $C_7$ by using the discharge method, another method for the list-coloring of planar graphs. Consequently, it may be concluded that planar graphs that are $C_k$-free for any $k\in \{3,4,5,6,7\}$ are 4-choosable. \\

\par In addition to these, the other conditions involved for being 4-choosable are the absence of intersecting triangles \cite{wang2002choosability}, absence of $C_4$ adjacent to $C_3$ \cite{borodin2008planar}, or every $C_5$ is not simultaneously adjacent to $C_3$ or $C_4$ \cite{xu2017sufficient}. Although the proof for 4-choosability of planar graphs without intersecting triangles uses both the boundary method and the discharge method, the result for planar graphs with $C_5$ not simultaneously adjacent to $C_3$ or $C_4$ is restricted to the discharge method. \\

\par The study of 3-choosability for planar graphs is subtler, since planar graphs are not 3-choosable in general.  In \cite{erdos1979choosability}, Erd\H{o}s et al. questioned the existence of a non-3-choosable planar bipartite graph. Although they identified a non-3-choosable planar graph, it was not bipartite. Alon \cite{alon1992colorings} answered this question by proving that every planar bipartite graph is 3-choosable which was a consequence of the finding that Alon-Tarsi number of such graphs is 3. The example of a non-3-choosable planar graph without 4-cycles, 5-cycles, and intersecting triangles appears in \cite{montassier2006bordeaux}. \\

\par The girth of a planar graph plays an important role when it comes to list-coloring. Thomassen \cite{thomassen19953} proved that every planar graph of girth at least 5 is 3-choosable. This can be seen as a list-coloring analogue, with a slightly stronger girth requirement,  of Gr{\"o}tzsch's theorem, which states that every planar graph of girth at least 4 is 3-colorable. Li \cite{li20093} extended Thomassen's result for planar graphs of girth at least 4 that has no $C_4$ that
shares a vertex with another $C_4$ or $C_5$. Later, Guo and Wang \cite{guo20113} proved a result that strengthened Li's statement. They showed that planar graphs of girth 4 that have no $C_4$ that shares a vertex with another $C_4$ or $C_5$ are 3-choosable.  \\

\par It is known that planar graphs with girth 4 are 4-choosable because they are 3-degenerate. Voigt \cite{voigt1995not} and Gutner \cite{gutner1996complexity} constructed planar graphs of girth 4 that are not 3-choosable, with 166 vertices and 164 vertices, respectively. For a planar graph of girth 4, 3-choosability is also guaranteed when it is either $\{C_5, C_6\}$-free \cite{lam20053} or $\{C_7, C_8\}$-free \cite{dvovrak2009planar}. \\

\par When there are no girth constraints, the results depend on forbidden cycles and often on the distance between the triangles. The distance between two triangles, say $T$ and $T'$, $d(T,T')=\min \{d(u,v) : u \in T, v \in T'\}$. Montassier et al. \cite{montassier2006bordeaux}  proved that a $\{C_4, C_5\}$-free planar graph is 3-choosable if $d(T,T')\geq 4$ for any $T,T'$ or is $C_6$-free with $d(T,T') \geq 3$ for any $T,T'$. Zhang and Sun \cite{zhang20083} proved that a $\{C_5, C_6\}$-free planar graph is 3-choosable if either it is $C_7$-free with $d(T,T')\geq 3$ for any $T,T'$ or $C_8$-free with $d(T,T')\geq 2$ for any $T,T'$. Wang et al. \cite{wang2008note} summarized the results in \cite{zhang2004three,zhang2005note,shen2007sufficient} and generalized it by stating that any planar graph without $C_4, C_i, C_j$ and $C_9$ where $i<j$ and $i,j \in\{5,6,7,8\}$ is 3-choosable.  \\
 
\par Apart from these, Grytczuk and Zhu \cite{grytczuk2020alon} showed that every planar graph contains a matching $M$ such that $G-M$ is 4-choosable as the Alon-Tarsi number of $G-M$ is at most 4. And Zhu \cite{zhu2023list} showed that if $L$ is a 4-list assignment on a planar graph with $|L(x) \cap L(y)|\leq 2$ for every edge $xy$, then $G$ is $L$-colorable.

\section{Toroidal graphs}\label{Toroidal graphs}
\hfill \par Toroidal graphs, the graphs that can be embedded on a torus, exhibit properties similar to planar graphs when it comes to coloring. By the Heawood bound for the torus, every toroidal graph is 7-colorable, and this bound is sharp since $K_7$ embeds on the torus. Toroidal graphs, being 6-degenerate, are 7-choosable. B{\"o}hme et al. \cite{bohme1999dirac} gave the necessary and sufficient condition for a toroidal graph to have list chromatic number 7. For a toroidal graph $G$, $\chi_\ell(G)=7$ if and only if $K_7$ is an induced subgraph of $G$. \\

\par Cai et al. \cite{cai2010choosability} gave the list-coloring of the toroidal graphs excluding cycles of fixed length $k$. Their results showed that a toroidal graph is 4-choosable if $k=3,4,5$, 5-choosable if $k=6$ and 6-choosable if $k=7$. Furthermore, they also proved that $\chi_\ell(G)=6$ if and only if $K_6$ is an induced subgraph of $G$. Based on this finding and the result given by B{\"o}hme et al.\cite{bohme1999dirac}, they conjectured that $\chi_\ell(G)=5$ if and only if $K_5$ is an induced subgraph of $G$. Choi \cite{choi2017toroidal} disproved this conjecture by giving an infinite family of graphs containing neither a $K_5$ nor a $C_6$ which is not 4-colorable, hence not 4-choosable. Alternatively, Choi \cite{choi2017toroidal} proved that a toroidal graph containing neither a $K_5-e$ nor a $C_6$ is 4-choosable. \\

\par While most of the results for planar graphs cannot be extended to toroidal graphs, there are some exceptions. The 4-choosability of toroidal graphs without intersecting triangles \cite{xiaofang20074} and the 4-choosability of toroidal graphs for which every $C_5$ is not adjacent simultaneously to $C_3$ or $C_4$ \cite{lv2018new} are two of such results. \\

\par The Cartesian product of two paths $P_m$ and $P_n$ is often called a grid and is denoted by $L_{m,n}$. Similarly, the Cartesian product of two cycles $C_m,C_n$ is often called a toroidal grid and is denoted by $T_{m,n}$, as it is a grid that can be embedded on a torus. Li et al. \cite{li2023alon} proved that for any $m,n \in \mathbb{N}$ the toroidal grid $T_{m,n}$ is 4-choosable. This provides another natural link between toroidal graphs and graph products.

\section{Perfect graphs}\label{Perfect graphs}
\hfill \par 
A graph $G$ is perfect if $\chi(H)=\omega(H)$ for every induced subgraph $H$ of $G$. Perfect graphs form a broad and structurally rich class. Among the earliest list-coloring results relevant to perfect graphs is the characterization of 2-choosable graphs by Erd\H{o}s, Rubin, and Taylor \cite{erdos1979choosability}, which implies that even cycles are chromatic-choosable. The same paper also treats odd-cycles, showing that every cycle is chromatic-choosable. \\

\par The next chromatic-choosable class of perfect graphs was chordal graphs \cite{tuza1996conjecture}. Chordal graphs are characterized by simplicial elimination ordering i.e., for any chordal graph with $n$ vertices, there exists a vertex ordering $v_1,v_2,\dots,v_n$ such that $N_{G_i}(v)$ induces a clique in $G_i$, where $G_i$ is the subgraph induced by $v_1,v_2,\dots,v_i$. The chromatic-choosability was immediate from this vertex ordering of the chordal graphs. \\

\par Hall's theorem on systems of distinct representatives (SDR) is one of the oldest and most common techniques used in list-coloring. An SDR for a family of sets $S_1, S_2, \dots, S_k$ is a sequence  $(x_1,x_2, \dots, x_k)$ such that $x_i \in S_i$ for each $i$, and  $x_i \neq x_j$ whenever $i\neq j$. We may treat the lists assigned for the vertices as the family of sets and a proper list-coloring as an SDR keeping the adjacency constraints. Using this idea, $\{\overline{K_3}, P_4\}$-free graphs have been proved to be chromatic-choosable \cite{woodall2001list}. \\

\par Interval graphs are yet another important class of perfect graphs. As a consequence of the Alon-Tarsi method, if the given graph $G$ has an acyclic orientation and a list assignment $L$ satisfying the condition that the outdegree of each vertex is at most one less than the size of list assigned to it, then $G$ is $L$-colorable. Woodall \cite{woodall2001list} observed that any interval graph $G$ has an acyclic orientation with $d^+(v) \leq \omega(G)-1$. Hence, $\chi_\ell(G)=\chi(G)=\omega(G)$. \\

\par An alternate proof for some classes of graphs like block graphs, $k$-trees and threshold graphs can be found in \cite{vasudevan2025chromatic}. Several important subclasses of perfect graphs --- including bipartite graphs, complete multipartite graphs, and claw-free perfect graphs --- have developed substantial list-coloring literatures. These classes are therefore treated separately in the following sections.

\section{Bipartite graphs}\label{Bipartite graphs}
\hfill \par This section focuses primarily on complete bipartite graphs, for which many exact and asymptotic results are known. Graphs with chromatic number at most two are precisely bipartite graphs; among graphs with at least one edge, the 2-chromatic graphs are exactly the bipartite graphs. This does not guarantee 2-choosability. The characterization for being 2-choosable uses the concept of $m$-core of a graph $G$, denoted by $\hat{G_m}$. It is the graph obtained by deleting the vertices in $G$ of degree less than $m$ successively until there are no more vertices of degree less than $m$. Specifically, a graph $G$ is 2-choosable if and only if $\hat{G_2}$ is of the form $K_1, C_{2m+2}$ or $\theta_{2,2,2m}$ for $m \geq 1$, where $\theta_{2,2,2m}$ is a graph consisting of two vertices $x$ and $y$ with three vertex disjoint $x,y$-paths each of length $2,2$ and $2m$ \cite {erdos1979choosability}. Recently, a generalization of this result was established in \cite{balachandran2021choice} that a graph $G$ is $k$-choosable if and only if $\hat{G_k}$ is $k$-choosable. \\

\par We first consider balanced complete bipartite graphs $K_{m,m}$. Erd\H{o}s et al. \cite{erdos1979choosability} showed that if $m=\binom{2k-1}{k}$ then $K_{m,m}$ is not $k$-choosable. Assign to the vertices in each part all the $k$ element subsets of $\{1,2,\dots,2k-1\}$ using each such subset exactly once in each part. Similarly for the other partite set, too. It can be seen that a list-coloring for this list assignment cannot be proper. They established $\chi_\ell(K_{m,m}) \leq \log_2m+3$. Using this result, they showed that the list chromatic number of any bipartite graph with equal sized partite sets, say $m$, is bounded below and above by $\log_6m$ and $3\log_6m$ respectively. They showed that asymptotically $\chi_\ell(K_{m,m}) \leq \log_2m+o(\log m)$. In their book on Combinatorial Nullstellensatz, Zhu and Balakrishnan
\cite{zhu2021combinatorial} gave a slightly better asymptote as $\log_2m+2$.

\subsection{$K_{1,q}$ and $K_{2,q}$}
\hfill \par Using the characterization of 2-choosable graphs, it is evident that the only 2-choosable graphs among the complete bipartite graph family are $K_{1,n}$ for $n \in \mathbb{N}$, $K_{2,2}$ and $K_{2,3}$, and all others fail to be chromatic-choosable. In their fundamental paper on list-coloring, Erd\H{o}s et al. \cite {erdos1979choosability} showed that any $K_{p,q}$ with $p \leq q$ is $(p+1)$-choosable and precisely $p$-choosable if and only if $q < p^p$. An immediate consequence of these results gives $\chi_\ell(K_{2,q})=3$ for $q>3$. 

\subsection{$K_{3,q}$}
\hfill \par Another immediate consequence of the result in \cite {erdos1979choosability} gives $\chi_\ell(K_{3,q})=3$ for $q<27$. Mahadev et al. \cite{mahadev19913}  extensively investigated the 3-choosability of bipartite graphs. They relied solely on the concept of the System of Distinct Representatives (SDR) to derive the results. They showed that $K_{3,q}$ for $q \geq27$ is not 3-choosable. But incorporating a later result that, $\chi_\ell(K_{p,q})= p+1$ when $q\geq p^p$ \cite{kaul2019list}, we can conclude that $\chi_\ell(K_{3,q})= 4$ for $q\geq27$.

\subsection{$K_{4,q}$}
\hfill \par The other details in \cite{mahadev19913} include the set of 3-choosable and non-3-choosable complete bipartite graphs. For $K_{4,q}$, there is a boundary value of 18 for q, less than or equal to which the graph becomes 3-choosable and greater than which the graph is no longer 3-choosable. A consequence of the result given by Hoffman and Johnson was used to investigate the result for a larger range of $q$. As $\chi_\ell(K_{p,q})= p$ when $(p-1)^{p-1}-(p-2)^{p-1} <q< p^p$ \cite{hoffman1993choice}, $\chi_\ell(K_{p,q})=4$ when $19 \leq q\leq 255$. 

\subsection{$K_{5,q}$}
\hfill \par  Mahadev et al. \cite{mahadev19913} showed that $K_{5,q}$ is 3-choosable for $q \leq 8$ and is not 3-choosable for $q \geq 15$,  leaving the cases $9 \leq q \leq 14$  unresolved. Shende and Tesman \cite{shende19953} partially filled the gap by showing that $K_{5,q}$ is 3-choosable for $9 \leq q \leq 12$.   West’s open problems \#15 cited a result of F{\"u}redi  that $K_{5,13}$ is not 3-choosable, and highlighted that its exact value was still open.   The case was resolved in \cite{mahadev19913} by showing that $\chi_\ell(K_{5,13})=4$ as cited by Hoffman and Johnson \cite{hoffman1996thwart}. 

\subsection{$K_{6,q}$}
\hfill \par Additionally, in \cite{mahadev19913} they added that $K_{6,q}$ for $q\geq 11$ is not 3-choosable. Although the 3-list-coloring for $q<11$ is not known, combining it with the result given by Hanson et al. \cite{hanson1996choosability} that $\chi_\ell(G)\leq 3$ for every bipartite graph with at most 13 vertices, we have $\chi_\ell(K_{6,6})=\chi_\ell(K_{6,7})=3$. In his preprint, O'Donnell established that $\chi_\ell(K_{6,11})=4$.\\

\par Hoffman and Johnson introduced a concept called thwart number in their paper \cite{hoffman1994restricted}. The thwart number of a graph $G$, $thw(G)$ is the smallest $k$ for which $G$ is not list-colorable when the lists assigned are $k$ subsets of $\{1,2,\dots,\chi_\ell(G)-1\}$. Subsequently, for any graph $G$, $\chi_\ell(G) \leq thw(G)$. In \cite{hoffman1996thwart}, the same authors obtained that $thw(K_{6,q})=3$ for $q=8,9,10$ and hence their list chromatic number of these graphs is also 3.

\subsection{$K_{7,q}$}
\hfill \par In \cite{mahadev19913}, they also established that $K_{7,q}$ is not 3-choosable. The result in  \cite{erdos1979choosability} given by $\chi_\ell(K_{7,7}) \geq 4$ supports this finding. The non-3-choosability of $K_{7,7}$ was given using the Fano plane configuration. The equality was established later. \\

\par We summarize the results into a table (refer Table \ref{list-coloring of Kpq}) for the ease of understanding.

\begin{table}[h!]
\centering
\begin{tabular}{|c|c|c|}
    \hline
    $K_{p,q}, q\geq p$ & $\chi_\ell(K_{p,q})$ & Value of $q$ \\
    \hline
    $K_{1,q}$ & 2 & $q \geq 1$ \\
    \hline
    $K_{2,q}$ & 2 & $q \leq 3$ \\
     & 3 & $q\geq 4$ \\
    \hline
    $K_{3,q}$ & 3 & $q \leq 26$\\
     & 4 & $q \geq 27$ \\
    \hline
    $K_{4,q}$ & 3 & $q \leq 18$ \\
     & 4 & $19 \leq q \leq 255$ \\
     & 5 & $q \geq 256$ \\
    \hline
    $K_{5,q}$ & 3 & $q \leq 12$\\
     & 4 & $q=13$\\
     & 5 & $176 \leq q \leq 3124$\\
     & 6 & $q \geq 3125$\\
    \hline
    $K_{6,q}$ & 3 & $q \leq 10$\\
     & 4 & $q=11$ \\
     & 6 & $2102 \leq q \leq 46655$\\
     & 7 & $q \geq 46656$\\
    \hline
    $K_{7,q}$ & 4 & $q=7$ \\
    & 7 & $31032 \leq q \leq 823542$\\
    & 8 & $q \geq 823543$\\
    \hline
\end{tabular}
    \caption{List chromatic numbers of $K_{p,q}$}
    \label{list-coloring of Kpq}
\end{table}

If $K_{p,q}$ is not $k$-choosable, then any $K_{p',q'}$ with $p'\geq p$ and $q'\geq q$ also fails to be $k$-choosable. Consequently, $\chi_\ell(K_{5,q})$ for $14 \leq q \leq 175$ is 4 or 5. Furthermore, the list-coloring problem for $K_{6,q}$ and $K_{7,q}$ remains unsolved $12 \leq q \leq 2101$ and $8 \leq q \leq 31031$, respectively, although both have a lower bound of 4 for the list chromatic number. For $p>7$, the list-coloring of $K_{p,q}$ is open for all $q\leq (p-1)^{p-1}-(p-2)^{p-1}$. We conclude the section with a conjecture on bipartite graphs.

\begin{conjecture}[Alon from Open Problem Garden] For a bipartite graph $G$, there exists a constant $c$ such that $\chi_\ell(G) \leq c \log \Delta(G)$.
\end{conjecture} 

\section{Multipartite graphs}\label{Multipartite graphs}
\hfill \par Currently, research on multipartite graphs is dominated by complete multipartite graphs. Throughout this section, $K_{r*n}$ denotes the complete multipartite graph with $n$ parts, each of size $r$. More generally, $K_{r_1*n_1,r_2*n_2, \dots}$  denotes the complete multipartite graph with $n_i$  parts of size $r_i$ for each $i$. In view of the conjecture posed by Ohba \cite{ohba2002chromatic}, He et al. \cite{he2008choice} gained interest in the chromatic-choosability of complete multipartite graphs which they posed as a conjecture later. We start with Ohba's conjecture.

\begin{conjecture}[\cite{ohba2002chromatic}] Any graph $G$ with a maximum of $2 \chi(G)+1$ vertices is chromatic-choosable.
\end{conjecture} 

\par Some weaker conditions of this conjecture have been proved. In the  same paper Ohba proved  chromatic-choosability of graphs whose order does not exceed $\chi(G)+ \sqrt{2\chi(G)}$ \cite{ohba2002chromatic}. Reed and Sudakov \cite{reed2003list}  proved two results in support of this conjecture. Initially, they showed that for any $0 < \varepsilon <1$, there exists a $n_0$ such that if $n_0 \leq n(G) \leq (2-\epsilon)\chi(G)$, then $G$ is chromatic-choosable. Also, they gave the chromatic-choosability of graphs with at most $\frac{5}{3} \chi(G)- \frac{4}{3}$ vertices. Ohba \cite{ohba2004choice} showed the equivalence of list chromatic number and chromatic number for any graph $G$ with at most $2 \chi(G)$ vertices and $\alpha(G) \leq 3$. Shen et al. \cite{shen2009ohba} improved this result under the same bound on the independence number for $n(G) \leq 2 \chi(G)+1$. Kostochka et al. \cite{kostochka2011ohba} refined this result further and proved the conjecture true for graphs with at most $2 \chi(G)$ vertices and $\alpha(G) \leq 5$.  \\

\par This remarkable conjecture was settled by Noel et al. \cite{noel2015proof}.  
They use the idea of minimal counterexamples and reformulate the coloring problem as a matching problem in a bipartite graph between vertices and colors.  Using a variant of Hall's marriage theorem, they analyze how colors can be distributed among parts of the complete multipartite graph. They show that any supposed obstruction must have extremely rigid structural properties which lead to a contradiction.
Zhu et al.\cite{zhu2022minimum} established that this bound is tight by putting forward the examples of $K_{3*(\frac{k}{2}+1),1*(\frac{k}{2}-1)}$ and $K_{4,2*(k-1)}$, both for even $k$, for the graphs with $2\chi(G)+2$ vertices and are non-chromatic-choosable. Additionally, they also showed that all other graphs with $2\chi(G)+2$ vertices are chromatic-choosable. \\

\par As any $k$-chromatic graph is the subgraph of some complete $k$-partite graph, He et al.\cite{he2008choice} gave another version of this conjecture as follows :

\begin{conjecture}[\cite{he2008choice}]
Any complete $k$-partite graph $G$ with $n(G)=2k+1$ is chromatic-choosable.
\end{conjecture} 

\subsection{$K_{r*n}$}
\hfill \par The first natural case in this notation is the complete multipartite graph with $n$ parts of size two.. This was given by Erd\H{o}s Rubin and Taylor  in their paper  the chromatic-choosability of $K_{2*n}$ \cite{erdos1979choosability} using  Hall's theorem. Later, Alon \cite{alon1992choice} used the probabilistic method and showed that $\exists c_1,c_2 \in \mathbb{N}$ such that $c_1n \log r \leq \chi_\ell(K_{r*n}) \leq c_2n \log r$. Gazit and Krivelevich \cite{gazit2006asymptotic} refined his work and gave more precise upper bounds for $\chi_\ell(K_{r*n})$ as $(1+o(1)) \frac{\log r}{\log(1+1/n)}$. In \cite{kierstead2000choosability}, Kierstead showed that $K_{3*n}$ is non-chromatic-choosable with list chromatic number $\lceil\frac{4n-1}{3}\rceil$. While Yang \cite{yang2003extension} showed that $\chi_\ell(K_{4*n})$ is bounded below and above by $\lfloor \frac{3n}{2}\rfloor$ and $\lceil \frac{7n}{4}\rceil$ respectively, Noel et al. \cite{noel2015beyond} improved the upper bound to $\lceil \frac{5n-1}{3}\rceil$ and Kierstead et al.\cite{kierstead2016choice} established that $\chi_\ell(K_{4*n})$ is equal to the lower bound given by Yang. Noel et al. \cite{noel2015beyond} showed that for smaller values of $k$, $K_{k*n}$ is bounded by $\lfloor \frac{2n(k-1)}{k}\rfloor$ and $\lceil \frac{n(k+1)-1}{3}\rceil$. Consequently, they showed that $\lfloor \frac{8n}{5}\rfloor \leq \chi_\ell(K_{5*n}) \leq 2n$ and $\lfloor \frac{5n}{3}\rfloor \leq \chi_\ell(K_{6*n}) \leq  \lceil  \frac{7n-1}{3}\rceil$.\\

\par Noel et al. \cite{noel2015beyond} showed that among the 3-chromatic graphs, $K_{3*k}$ has the largest list chromatic number. Based on this result, they conjectured that $K_{m*k}$ has the maximum list chromatic number among all the $k$-chromatic graphs.

\begin{conjecture}[Noel from Open Problem Garden] If $G$ is a $k$-chromatic graph with at most $mk$ vertices, then $\chi_\ell(G) \leq \chi_\ell(K_{m*k})$.
\end{conjecture}

\subsection{$K_{m*1,2*(n-1)}$}
\hfill \par Gravier and Maffray \cite{gravier1998graphs} proved that any $k$ colorable with at most two color classes having 3 vertices and all other color classes having at most 3 vertices is $k$-choosable. Using this result, it can easily be deduced that $\chi_\ell(K_{3*1,2*(n-1)})=n$. Enomoto et al. \cite{enomoto2002choice} proved the results for $m=4,5$. They showed that 
$\chi_\ell(K_{4*1,2*(n-1)})$ is $n+1$, if $n$ is even while it is $n$, if $n$ is odd, and $\chi_\ell(K_{5*1,2*(n-1)})=n+1$. Later, Vetr{\'\i}k \cite{vetrik2012list} showed that $K_{\frac{(t+2)(t+3)}{2}*1,2*(n-1)}$ is $(n+t)$-choosable. By putting $t=1$, and as $\chi_\ell(K_{5*1,2*(n-1)})=n+1$, we can conclude that $\chi_\ell(K_{6*1,2*(n-1)})=n+1$. They also established that
$K_{\binom{2k}{k}^2*1,2*(n-1)}$ is tightly bounded below by $n+k$ whenever $n=pk+1$ for some odd integer $p$.

\subsection{$K_{r_1*n_1,r_2*n_2, \dots}$}
\hfill \par Using the same idea that used to show the chromatic-choosability of $K_{3*1,2*(n-1)}$, it is immediate that $K_{3*2,2*(n-2)}$, $K_{3*2,2*(n-k-2),1*k}$ and $K_{3*1,2*(n-k-1),1*k}$ are chromatic-choosable. Later,  Ohba \cite{ohba2004choice} in his paper generalized this result for all possible combinations of partite sets of size 3, 2 and 1 and hence proved the chromatic-choosability of $K_{3*k_1,2*(n-k_1-k_2),1*k_2}$. Along with it, he also established that $\chi_\ell(K_{3*k_1, 1*k_2})= \max \{k_1+k_2, \lceil \frac{4k_1+2k_2-1}{3} \rceil\}$. Through close observation, it can be easily seen that the terms whose maximum is to be compared are nothing but $\chi(G)$ and $\lceil \frac{n(G)+\chi(G)-1}{3} \rceil\}$. Based on this result, Noel et al. \cite{noel2015beyond} gave the upper bound for the list chromatic number that we found in Section 2. Shen et  al.\cite{shen2008choosability} established the chromatic-choosability of $K_{4*1,3*1,2*(n-4),1*2}$, $K_{5*1,3*1,2*(n-5),1*3}$ and $K_{(k+3)*1,2*(n-k-1),1*k}$ in support of Ohba conjecture. 

\section{Claw-free graphs}\label{Claw-free graphs}
\hfill \par \textit{Claw-free} graphs are $K_{1,3}$-free graphs. Gravier and Maffray \cite{gravier1997choice} initiated the systematic study of chromatic-choosability for claw-free graphs and formulated the following conjecture.

\begin{conjecture}[LCC for claw-free graphs \cite{gravier1997choice}]
For every claw-free graph $G$, $\chi_\ell(G)=\chi(G)$.
\end{conjecture}

\par Gravier and Maffray \cite{gravier1997choice} proved that any $k$-colorable claw-free graph with at most $2k+2$ vertices is $k$-choosable \cite{gravier1998graphs}. Chudnovsky and Seymour \cite{chudnovsky2010claw} gave a slightly better bound for the list chromatic number of claw-free graph as the maximum among $\chi(G)$ and $\frac{n(G)}{2}$. They showed that for claw-free graphs $G$ with $\chi_\ell(G-\{v\}) \leq 2\omega(G-\{v\})$, for some vertex $v$ in $G$ whose neighborhood is the union of two cliques, then $\chi_\ell(G)\leq 2\omega(G)$. In the same paper, the authors showed that twice the clique number bounds the chromatic number of a claw-free graph with independence number at least three. As a consequence of this bound, a \{claw, $K_4$\}-free graph is 6-colorable and a \{claw, $K_5$\}-free graph is 8-colorable. 
Esperet et al. \cite{esperet2014list} studied the corresponding list-coloring questions and proved that \{claw, $K_4$\}-free graphs are 4-choosable and \{claw, $K_5$\}-free graphs are 7-choosable.
\\

\par It is well known that any graph is $(\Delta(G)+1)$-colorable. Borodin and Kostochka \cite{borodin1977upper} had conjectured that every graph $G$ with $\Delta(G) \geq 9$ and that has no clique of size $\Delta(G)$ is $(\Delta(G)-1)$-colorable. Cranston and Rabern \cite{cranston2017clawfree} proved that the conjecture is true for claw-free graphs. In the same paper by these authors, we can find the list-coloring analogue of the conjecture given by Borodin and Kostochka.

\begin{conjecture}[\cite{cranston2017clawfree}]
Every graph $G$ with $\Delta(G) \geq 9$ and has no clique of size $\Delta(G)$ is $(\Delta(G)-1)$-choosable.
\end{conjecture}

\par In a later paper, the same authors have shown that if $G$ is a claw-free graph with $\Delta(G) \geq 69$ and has no clique of size $\Delta(G)$ is $(\Delta(G)-1)$-choosable. For the graphs for which $9 \leq \Delta(G) \leq 68$ and the same condition, the conjecture is still open. \\

\par Although chromatic-choosability remains largely open for claw-free graphs in general, stronger results are known for claw-free perfect graphs.

\subsection{Claw-free perfect graphs}

\hfill \par We first discuss elementary graphs, a subclass of claw-free perfect graphs. A graph $G$ is said to be an \textit{elementary graph} if its edges can be bicolored such that every induced $P_3$ is bicolored. Gravier and Maffray \cite{gravier1997choice} proved that any elementary graph $G$ with clique number at most 3 is chromatic-choosable. \\ 

\par Another subclass of claw-free perfect graphs is peculiar graphs. Rather than a definition, we have a construction for \textit{peculiar graph} which is as follows. Let $K_n$ be a complete graph whose set of vertices is split into six pairwise disjoint nonempty sets $A_1, B_1, A_2, B_2, A_3, B_3$. For each $i= 1, 2, 3,$ remove at least one edge with one endpoint in $A_i$ and the other endpoint in $B_{i+1}$ (subscript is modulo 3). Finally, add three pairwise disjoint nonempty cliques $Q_1,Q_2,Q_3$, and, for each $i= 1, 2, 3$, provide adjacencies between every vertex in $Q_i$ and $K_n-(A_i \cup B_i)$. Gravier and Maffray \cite{gravier1998graphs} established the chromatic-choosability of peculiar graphs.  \\

\par The classification of claw-free perfect graphs according to Chv{\'a}tal \cite{chvatal1988recognizing} states that every claw-free perfect graph that has no clique cutset is either elementary or peculiar. A graph $G$ is said to have a \textit{clique cutset} if it has a vertex cut that induces a clique. Giving an extension to 3-chromatic claw-free perfect graphs with clique cutset, Gravier and Maffray \cite{gravier2004choice} established that every 3-chromatic claw-free perfect graph is 3-choosable. \\

\par After nearly two decades, Gravier et al. \cite{gravier2016choosability} contributed further to the chromatic-choosability of claw-free perfect graphs by proving the same for claw-free perfect graphs with clique number 4. The authors formally stated the chromatic-choosability of claw-free perfect graphs as a conjecture in the same paper before giving the proof for their result.

\begin{conjecture}[LCC for claw-free perfect graphs \cite{gravier2016choosability}]
For every claw-free perfect graph $G$, $\chi_\ell(G)=\chi(G)$.
\end{conjecture}

\par The conjecture for claw-free perfect graphs is a natural intermediate problem toward the broader list-coloring conjecture for claw-free graphs. Using the classification of claw-free perfect graphs, we can say that the conjecture remains open for non-peculiar claw-free perfect graphs with clique number at least five.

\subsection{Line graphs}
\hfill \par We considered different versions of the list-coloring conjecture. We now  discuss two more such conjectures. The first one was stated independently by Vizing, by Gupta, and by Albertson and Harris. The second is its formulation for line graphs by H{\"a}ggkvist and Chetwynd \cite{haggkvist1992some}. \\

\par The \textit{edge coloring} of a graph is the coloring of its edges in such a way that no two incident edges receive the same color. We say that a given graph is \textit{$k$-edge colorable} if the graph admits an edge coloring using $k$ colors. The least $k$ for which a graph $G$ is $k$-edge colorable is called the \textit{chromatic index} of $G$, denoted by $\chi'(G)$.\\

\par Let $L$ be a list assignment on the edge set of $G$. If $G$ can be edge colored from their corresponding lists, we say that $G$ is \textit{$L$-edge colorable}. If all lists are of the same size, say $k$ and $G$ is $L$-edge colorable for every $k$-list assignment $L$ on $G$, then $G$ is \textit{$k$-edge choosable}. The least $k$ for which a graph $G$ is $k$-edge choosable is called the \textit{list chromatic index} of $G$, denoted by $\chi_\ell'(G)$ or $ch'(G)$.

\begin{conjecture}[List Edge Coloring Conjecture (LECC)] For any graph $G$, $\chi_\ell'(G)=\chi'(G)$.
\end{conjecture}

Corresponding to any given graph $G$, there is a graph associated with it, called the line graph of $G$. The \textit{line graph} of $G$, denoted by $L(G)$ is the graph with $V(L(G))=E(G)$ and any pair of vertices in $L(G)$ are adjacent if the corresponding pair of edges share a common endpoint. Line graphs form a  subclass of claw-free graphs. \\

\par Observe that the edge coloring of a graph is equivalent to the vertex coloring of its line graph. Similarly, the list edge coloring of a graph is equivalent to the list-coloring of its line graph. Thus, H{\"a}ggkvist and Chetwynd rephrased the LECC as follows. 

\begin{conjecture}[LCC for Line Graphs \cite{haggkvist1992some}] 
For any graph $G$, $\chi_\ell(L(G))=\chi(L(G))$.
\end{conjecture}

\par Several results have been proved in favour of the conjecture. Snarks \cite{harris1985problems}, $K_{3,3}$ \cite{chetwynd1989note}, $K_{4,4}$ \cite{alon1992colorings}, $K_{6,6}$ \cite{alon1992colorings}, bipartite graphs \cite{galvin1995list}, 2-connected regular planar graphs \cite{ellingham1996list}, $d$-regular multiplanar graphs \cite{ellingham1996list}, $K_{2n+1}$ \cite{haggkvist1997new}, line-perfect graphs \cite{peterson1999edge}, series-parallel graphs \cite{juvan1999list}, multicircuits \cite{woodall1999edge}, multigraphs with a vertex whose deletion makes the graph bipartite \cite{plantholt1999list},  graphs with no circuit of size 4 or more \cite{gutner2009some}. \\

 We highlight Galvin's proof for the bipartite graphs  \cite{galvin1995list} as it introduced one of the most important applications of the kernel method to list coloring. A {\em kernel } in a digraph $D$ is an independent set $S\subseteq V(D)$ such that every vertex in $V(D)\setminus S$ has an outgoing edge to a vertex of $S$. A digraph is {\em kernel-perfect } if every induced subdigraph has a kernel. The kernel lemma states that if a graph $G$ has a kernel-perfect orientation $D$ satisfying $d_D^+(v)\leq |L(v)|-1$ for every vertex $v$, then $G$ is $L$-colorable. \\
 
 Let $H$ be a bipartite multigraph with bipartition $(X, Y)$, and let $G = L(H)$. The vertices of $G$ are the edges of $H$. For $x \in X$, the set of edges incident with $x$ may be viewed as a row, and for $y \in Y$, the set of edges incident with $y$ may be viewed as a column.  Starting from a proper $n$-edge-coloring $f: E(H) \to \{1, \ldots, n\}$, where $n=\Delta(H)$,  Galvin defines an orientation of $L(H)$ as follows: if two adjacent edges $u, v$ lie in the same row, orient $u \to v$ whenever $f(u) > f(v)$; if they lie in the same column, orient $u \to v$ whenever $f(u) < f(v)$. Thus, the orientation is determined by the edge-coloring together with the bipartition of $H$.\\

This orientation has maximum outdegree at most $n - 1$. Moreover, it is a normal orientation of the line graph, and Galvin uses the stable marriage theorem, in Maffray’s formulation, to show that every induced subdigraph has a kernel. The kernel lemma then implies that $L(H)$ is $n$-choosable. Since the chromatic number of $L(H)$ is the chromatic index of $H$, and a bipartite multigraph has chromatic index equal to its maximum degree, it follows that:
\[
\chi_\ell(L(H)) = \chi(L(H)).
\]
Equivalently, the list chromatic index of every bipartite multigraph equals its chromatic index.

\section{Square of graphs}\label{Square of graphs}
\hfill \par The \textit{square} of a given graph $G$, denoted by $G^2$ is the graph with $V(G^2)=V(G)$ and $E(G^2)$=\{$uv:u$ and $v$ are of distance at most 2 in $G$\}. This section surveys list-coloring results for squares of graphs. We begin with total coloring, since total coloring can be naturally reformulated as a coloring problem for squares of subdivision graphs. \\

\par The \textit{total coloring} of a graph is the coloring of its vertices and edges in such a way that no two adjacent vertices, no two incident edges and,  no incident vertex-edge pair receive the same color. We say that a given graph is $k$-total colorable if the graph admits a total coloring using $k$ colors. The least $k$ for which a graph $G$ is $k$-total colorable is called the \textit{total chromatic number} of $G$, denoted by $\chi''(G)$. \\

\par Let $L$ be a list assignment on the vertex set and edge set of $G$. If $G$ can be totally colored from their corresponding lists, we say that $G$ is \textit{$L$-totally colorable}. If all lists are of the same size, say $k$ and $G$ is $L$-totally colorable for every $k$-list assignment $L$ on $G$, then $G$ is \textit{$k$-totally choosable}. The least $k$ for which a graph $G$ is $k$-totally choosable is called the \textit{list total chromatic number} or \textit{total choosability} of $G$, denoted by $\chi_\ell''(G)$ or $ch''(G)$.\\

\par Borodin et al. \cite{borodin1997list} stated a conjecture on the coincidence of total chromatic number and list total chromatic number in 1996. This is called the list total coloring conjecture.

\begin{conjecture}[List Total Coloring Conjecture (LTCC) \cite{borodin1997list}]
For any multigraph $G$, $\chi_\ell''(G)=\chi''(G)$.
\end{conjecture}

\par Corresponding to any given graph $G$, there is a graph associated with it, called the total graph of $G$. The \textit{total graph} of $G$, denoted by $T(G)$ is the graph with $V(T(G))=V(G) \cup E(G)$ and any pair of vertices in $T(G)$ are adjacent if the corresponding pair of vertices are adjacent or the corresponding pair of edges are incident or the corresponding edge is incident on the corresponding vertex. Note that the total coloring of a graph is equivalent to the coloring of its total graph. Analogously, the list total coloring of a graph is equivalent to the list-coloring of its total graph. Thus, LTCC can be rephrased as $\chi_\ell(T(G))=\chi(T(G))$. \\

\par  We do not discuss LTCC in detail, but we recall several classes of graphs for which LTCC holds. Kostochka and Woodall showed that for some graphs like multigraphs whose underlying simple graph is a forest, multigraphs whose underlying simple graph is a circuit and bipartite graphs with total chromatic number $\Delta+2$, LTCC holds. See \cite{kostochka2001choosability, kostochka2002total1, kostochka2002total2} for details. While working with these graphs, they observed that for any multigraph $G$, $T(G)=H^2$ where $H$ is the graph obtained by subdividing every edge of $G$ exactly once. At this point, Kostochka and Woodall \cite{kostochka2001choosability} conjectured the chromatic-choosability of squares of graphs. 

\begin{conjecture}[List Square Coloring Conjecture (LSCC) \cite{kostochka2001choosability}] For any graph $G$, $\chi_\ell(G^2)=\chi(G^2)$.
\end{conjecture}

\par It can be seen that LSCC implies LTCC. Since $T(G)=H^2$, where H is obtained from $G$ by subdividing every edge exactly once, the List Total Coloring Conjecture is a special case of the List Square Coloring Conjecture on bipartite graphs of this type. Hence, LSCC is true for all graphs for which LTCC holds. A look-alike of this conjecture and a small variant of this question is raised in \cite{kim2014bipartite}. \\

\noindent \textbf{Question \cite{kim2014bipartite}:} \textit{If $G$ is a bipartite graph such that every vertex of one partite set has degree at most 2, then is it true that $\chi_\ell(G^2)= \chi(G^2)$?}\\

\noindent \textbf{Question \cite{kim2014bipartite}:} \textit{If the above question is true, then what is the largest $k$ such that $G^2$
is chromatic-choosable for every bipartite graph $G$ with a partite set in which each vertex has degree at most $k$?}\\

\par Later, the conjecture was disproved by Kim and Park \cite{kim2015counterexamples}. For each prime $p>2$, they constructed a graph $G$ with $2p^2-p$ vertices using Latin squares. They constructed the graph in such a way that $G^2$ is the multipartite graph $K_{p*(2p-1)}$ and $\chi_\ell(K_{p*(2p-1)})-\chi(K_{p*(2p-1)}) \geq p-1$. Later, Kim and Park constructed a bipartite graph that violates the LSCC \cite{kim2014bipartite}. For each prime $p \geq 7$, it was a graph with $p(p^2+2p-1)$ vertices and they proved that the $jump(G^2)>p^2-6p+3$, is arbitrarily large for larger values of $p$.\\

\par We look at two open problems posed by Noel on the upper bound for list chromatic number of square of a graph. \\

\noindent \textbf{Open Problem (Noel from Open Problem Garden):} \textit{Does there exist a function $f$ such that for every graph $G$, $\chi_\ell(G^2) \leq f(\chi(G^2))$?}\\

\noindent \textbf{Open Problem (Noel from Open Problem Garden):} \textit{Does there exist a constant $c$ such that $\chi_\ell(G^2) \leq c\chi(G^2) \log \chi(G^2)$?}

\subsection{Squares of planar graphs}
\hfill \par The squares of planar graphs have also been an area of extensive research. For squares of planar graphs, the maximum degree plays a central role in known upper bounds. Wegner \cite{wegner1977graphs} conjectured the following bounds for $\chi(G^2)$: 
\begin{center}
    $\chi(G^2) \leq 
    \begin{cases}
        7, \Delta \leq 3\\
        \Delta+5, 4 \leq \Delta \leq 7\\
        \lfloor \frac{3}{2} \Delta \rfloor+1, \Delta \geq 8.
    \end{cases}$
\end{center}
Havet et al. \cite{havet2007list} proved the third case asymptotically as $\chi_\ell(G^2) \leq \frac{3}{2} \Delta(G)[1+o(1)]$. They formally proposed the list-coloring conjecture for the square of planar graphs.

\begin{conjecture}[LSCC for planar graphs \cite{havet2007list}]
For any planar graph $G$, $\chi_\ell(G^2)=\chi(G^2)$.
\end{conjecture}

\par It was later disproved by Hasanvand \cite{hasanvand2022list}. He gave a set of infinite families of graphs whose squares are not chromatic-choosable. This includes bipartite cubic graphs with girth 6, planar bipartite cubic graphs, and claw-free planar cubic graphs. He also established that there exists a graph $G$ such that $jump(G^2)\geq n$, where $n$ is an arbitrary positive integer.\\

\par Let $G_{g=m, \Delta=k}$ denote a planar graph with girth $m$ and maximum degree $k$. If the subscripts have inequalities for $g$ and/or $\Delta$, then girth and maximum degree will have the respective conditions. Borodin et al.\cite{borodin2008list} gave certain bounds for girth and maximum degree under which the square of the graph has list chromatic number $\Delta+1$. The graphs are as follows : $G_{g=7, \Delta\geq 30}$, $G_{g\geq 8, \Delta\geq 15}$, $G_{g=10, \Delta\geq 9}$, $G_{g\geq 11, \Delta\geq 7}$, $G_{g\geq 12, \Delta=6}$, $G_{g\geq13, \Delta=5}$, $G_{g\geq15, \Delta=4}$ and $G_{g\geq 24, \Delta=3}$. They also showed that for girth less than 7, there exists a planar graph $G$ such that $\chi_\ell(G^2)=\Delta(G)+2$, for arbitrarily large $\Delta(G)$. Ivanova \cite{ivanova2011list} strengthened the conditions given by Borodin et al. \cite{borodin2008list}. Ivanova strengthened these conditions by proving $G_{g\geq 7, \Delta\geq 16}$, $G_{8 \leq g\leq 9, \Delta\geq 10}$, $G_{10\leq g\leq 11, \Delta\geq 6}$ and $G_{g\geq 12, \Delta=5}$ have  $\chi_\ell(G^2)=\Delta(G)+1$. \\

\par Borodin and Ivanova  had worked consistently on the $(\Delta+2)$-list-coloring of squares of planar graphs. They kept the bound on the girth constant and kept tightening the bound on the maximum degree. From $G_{g\geq 6, \Delta\geq 36}$, they improved it to $G_{g\geq 6, \Delta\geq 24}$ and then to $G_{g\geq 6, \Delta\geq 18}$. These works appear in \cite{borodin2009list, borodin20092, borodin20091}. Bonamy et al.\cite{bonamy2014graphs} tightened the bound on maximum degree further to 17 and showed that $\chi_\ell(G_{g\geq 6, \Delta\geq17}^2) \leq \Delta(G)+2$. Keeping the girth bound the same and tightening the maximum degree bound to 9 would make the square of the planar graph $(\Delta+3)$-list-colorable. This was the result given by Bu and Shang \cite{bu2016list} that $\chi_\ell(G_{g\geq 6, \Delta\geq9}^2) \leq \Delta(G)+3$. Cranston and Jaeger \cite{cranston2017planar} established that if $G$ is a $\{C_4, C_5\}$-free planar graph with $\Delta(G) \geq 32$, then $\chi_\ell(G^2) \leq \Delta(G)+3$.\\

\par Jin and Miao \cite{jin2022list} worked on planar graphs with girth at least 5 and showed that $\chi_\ell(G_{g\geq 5, \Delta\geq40}^2) \leq \Delta(G)+4$. Bu and Yan \cite{bu2015list} proved that $\chi_\ell(G_{g\geq 5, \Delta\leq5}^2) \leq 13$. They also added that planar graphs $G$ with $\Delta(G) \geq 12$ and without 3, 5-cycles and intersecting 4-cycles have $\chi_\ell(G^2) \leq \Delta(G)+6$; . Bu and Shang \cite{bu2016list} replaced the condition of forbidden cycles and intersecting cycles using the condition on girth and established that $\chi_\ell(G_{g\geq 5, \Delta\geq12}^2) \leq \Delta(G)+6$. \\

\par For every planar graph with neither 3-cycles nor intersecting 4-cycles and $\Delta(G) \geq 18$, $\chi_\ell(G^2) \leq \Delta(G) + 8$ \cite{zhou2022list}. Cranston et al. \cite{cranston2013choosability} established that $\chi_\ell(G_{g\geq 11, \Delta=4}^2) \leq 6$. Zhou and Sun \cite{zhou20212} proved that $\chi_\ell(G_{g\geq 10, \Delta\leq4}^2) \leq 7$ and $\chi_\ell(G_{g\geq 6, \Delta\leq4}^2) \leq10$. \\

\par The example of the 3-colorable non-4-choosable planar graph given by Mirzakhani \cite{mirzakhani1996small} satisfied the inequality $\chi_\ell(G) \geq \chi(G)+2$. In light of this inequality, Hasanvand questioned the existence of a planar graph whose square satisfies this inequality. \\

\noindent \textbf{Question \cite{kim2014bipartite}:} \textit{Is there a planar graph G satisfying $\chi_\ell(G^2) \geq \chi(G^2)+2$?}

\subsection{Squares of subcubic planar graphs}
\par Some cubic planar graphs provide counterexamples to the list square coloring conjecture for planar graphs. We now focus on results in this class. A graph with $\Delta \leq 3$ is called a \textit{subcubic graph}. Thomassen \cite{thomassen2018square} showed that the square of a subcubic planar graph is 7-colorable. This naturally raises the question whether the square of a subcubic planar graph is 7-choosable. As in the planar case, we can see that girth is an important factor for the results involving squares of subcubic planar graphs.\\

\par In \cite{cranston2008list}, Cranston and Kim established that for every subcubic graph $G$ except the Petersen graph, $G^2$ is 8-choosable. This bound is tight, as the square of the graph obtained by deleting an edge from the Petersen graph is $K_8$. They also proved that the square of every subcubic planar graph of girth at least 7 is 7-choosable, and that the square of every subcubic planar graph of girth at least 9 is 6-choosable. An improvement to this result was given by Kim and Lian in \cite{kim2024square} that the square of a subcubic planar graph of girth at least 6 is 7-choosable. Dvo{\v{r}}{\'a}k et al. \cite{dvovrak2008list} extended these results to the squares of a subcubic planar graph of girth at least 14 and at least 24 are 5-choosable and 4-choosable respectively. Havet \cite{havet2009choosability} made an improvement to one of the results given by Dvo{\v{r}}{\'a}k et al. and thus establishing that the square of a subcubic planar graph of girth at least 13 is 5-choosable. It can be observed that as girth of a subcubic planar graph $G$ is increased, the bound on $\chi_\ell(G^2)$ is tightened gradually. \\

\par Few recent works on list-coloring of squares of subcubic planar graphs revolve around forbidden cycles. Jin et al. \cite{jin2026square} showed that if $G$ is a $\{C_4, C_5\}$-free subcubic planar graph, $G^2$ is 7-choosable. While Kim et al. \cite{kim2025square} showed that the absence of  $C_4$is not a necessary for the square of the a subcubic planar graph to be 7-choosable. Kim and Luo \cite{kim2025squares} showed that $G^2$ is 6-choosable for any
$C_k$-free subcubic planar graphs for $4\leq k \leq 8$.

\section{Powers of graphs}\label{Powers of graphs}
\hfill \par The \textit{$k\textsuperscript{th}$ power} of a given graph $G$, denoted by $G^k$ is the graph with $V(G^k)=V(G)$ and $E(G^k)$=\{$uv:u$ and $v$ are of distance at most $k$ in $G$\}. Prowse and Woodall \cite{prowse2003choosability} were  among the first to study the list-coloring of powers of graphs. They showed that the powers of cycles are chromatic-choosable. \\

\par In view of the counterexample given by Kim and Park for LSCC, Zhu questioned whether there exists an integer $k$ such that every $k^{th}$ power of a graph is chromatic-choosable. In \cite{kim2015chromatic}, Kim et al. showed that no such $k$ exists. Moreover, for every odd $k$, there exists a graph $G$ such that $\chi_\ell(G^k)=\lceil \frac{4\chi(G^k)-1}{3} \rceil$. They established this by showing that for any $k \geq 2$, there exists a graph $G$ such that $\chi_\ell(G^k) \geq \frac{10}{9} \chi(G^k)-1$. In a related direction, Kosar et al. \cite{kosar2014note} showed that there exists a constant $c$ such that for any $k \in \mathbb{N}$, there exists an infinite family of graphs $G$ with $\chi_\ell(G^k)$ unbounded such that $\chi_\ell(G^k) \geq c\chi(G^k) \log \chi(G^k)$. They questioned the existence of a constant $c'$ such that for any $k \in \mathbb{N}$, there exists a graph $G$ such that $\chi_\ell(G^k) \leq c'\chi(G^k) \log \chi(G^k)$, which is a version of the question posed by Noel without any constraints on power. They have also given upper bounds for the list chromatic number of the power of a graph $G$ as, $\chi_\ell(G^k) \leq (\chi(G^k))^3$. In particular, if $k \geq 3$ and $k$ is odd, then $\chi_\ell(G^k) \leq \Delta(G)(\chi(G^k))^2$. \\

\par Kim et al. \cite{kim2015chromatic} posed two interesting questions towards the end of their paper. \\

\noindent \textbf{Question \cite{kim2015chromatic} :} \textit{If $G$ is claw-free, does there exist a constant $k$ such that $G^k$ is chromatic-choosable?}\\

\noindent \textbf{Question \cite{kim2015chromatic} :} \textit{If $G$ is claw-free, is $G^k$ chromatic-choosable for every $k \geq 2$?}\\

We conclude this section with the chromatic-choosability of powers of path graphs. Kaul and Mudrock \cite{kaul2018alon} had shown that $\chi_\ell(P_n^r)=\chi(P_n^r)= r+1$ for all $1 \leq r \leq n-1$.

\section{Graph operations}\label{Graphs operations}
\hfill \par This section surveys list-coloring results for standard graph operations, including joins, Cartesian products, and lexicographic products.

\subsection{Join of graphs}
\hfill \par The join of two graphs $G$ and $H$, denoted by $G \lor H$, is the graph with vertex set $V(G\lor H)=V(G) \cup V(H)$ and $E(G\lor H)=E(G) \cup E(H) \cup \{uv : u \in V(G), v \in V(H)\}$. It is known that $\chi(G \lor H)=\chi(G)+\chi(H)$. But this need not necessarily hold for the list chromatic number. If $S_k$ denotes a stable set of size $k$, we have $K_{n,n^n}=S_n \lor S_{n^n}$. When $\chi_\ell(S_n \lor S_{n^n})=n+1$, $\chi_\ell(S_n)+ \chi_\ell(S_{n^n})= 2$. Replacing $S_n$ with any graph $G$ with $n$ vertices, we can deduce that $G \lor S_{n^n}$ is not $n$-choosable.\\

\par Ohba \cite{ohba2002chromatic} gave a condition under which the join of two graphs is chromatic-choosable. He proved that whenever $n(G)+n(H) \leq (\chi(G)+\chi(H))+ 2\sqrt{(\chi(G)+\chi(H))}$, then $\chi_\ell(G \lor H)= \chi(G \lor H)$. He also
proved that for any graph $G$, there exists a non-negative integer $n_0$ such that for any integer $n\geq  n_0$, $G \lor K_n$ is chromatic-choosable. Moreover, $\chi_\ell(G \lor K_1)\leq \chi_\ell(G)+1$, and the equality need not hold always. \\

\par Gravier et al. \cite{gravier2003list} showed that for a graph $G$, there exists a function $f(G)$ defined  as the smallest integer $k$ such that $G \lor S_k$ is not $n(G)$-choosable. Using this result, it follows that $f(K_n)=1$ and $f(S_n)=n^n$. Thus, for any graph with $n$ vertices, $1 \leq f(G) \leq n^n$. They have also given the value of $f$ for some graphs, such as complete graphs, graphs that are disjoint unions of cliques, and triangle-free graphs.\\

\noindent \textbf{Open Problem (Allagan from Open Problem Garden):} \textit{Given $a,b\geq2$, what is the smallest $t\geq 0$ such that $\chi_\ell(K_{a,b}\lor K_t) \leq \chi(K_{a,b}\lor K_t)$?}

\subsection{Cartesian product}
\hfill \par The Cartesian product of two graphs $G=(V(G),E(G))$ and $H=(V(H),E(H))$, denoted by $G \Box H$, is the graph with vertex set $V(G) \times V(H)$ and any two vertices $(u,v),(u',v')$ are adjacent if and only if either $u=u'$ and $vv' \in E(H)$ or $v=v'$ and $uu' \in E(G)$. Sabidussi \cite{sabidussi1957graphs} proved that $\chi(G \Box H)=\max\{\chi(G),\chi(H)\}$. Hence one possible question was $\chi_\ell(G \Box H)=\max\{\chi_\ell(G),\chi_\ell(H)\}$, which was proved false by Borowiecki and Jozef \cite{borowiecki2006list}. \\

\par As $K_m \Box K_n \cong L(K_{m,n})$, it was trivial from Galvin's result that \cite{galvin1995list}, $\chi_\ell(K_m \Box K_n)=\chi(K_m \Box K_n)=\max\{m,n\}$. One of the first  contributions to list-coloring of products of graphs was given by Borowiecki and Jozef \cite{borowiecki2006list} in the form of an upper bound for the Cartesian product of two graphs as
    \begin{center}
        $\chi_\ell(G \Box H) \leq \min\{\chi_\ell(G) + $col$(H), $col$(G) + \chi_\ell(H)\}-1$
    \end{center}
and an immediate generalization for the Cartesian product of $k$ graphs as
    \begin{center}
        $\chi_\ell(G_1 \Box G_2 \Box \dots G_k) \leq \chi_\ell(G_1) + $col$(G_2)+\dots+$col$(G_k)-(k-1)$
    \end{center}
They have also shown that, for any graph $G$ with $k$ vertices, $\chi_\ell(G \Box K_{m,n}) = \chi_\ell(G) + m$, where $n=(m+\chi_\ell(G)-1)^{mk}$. And $\chi_\ell(G \Box G) = 2k+1$ when $G=K_{k,(2k)^{k(k+k^k)}}$. The authors have posed two conjectures as well.\\

\begin{conjecture}[\cite{borowiecki2006list}]
For every pair of graphs $G_1$ and $G_2$, there exists a constant $C$ such that $\chi_\ell(G_1 \Box G_2) \leq C(\chi_\ell(G_1)+\chi_\ell(G_2))$.
\end{conjecture}

\begin{conjecture}[\cite{borowiecki2006list}]
Let $G_1$ and $G_2$ be two graphs with $\Delta(G_1) \leq \Delta$ and $\Delta(G_2) \leq \Delta$. Then $\chi_\ell(G_1 \Box G_2) \leq \Delta+o(\Delta)$.
  \end{conjecture}
\par Kaul and Mudrock \cite{kaul2018alon} have given a few results for the Cartesian product of two graphs in which one of them is a path graph. They proved that the Cartesian product of a path and a cycle is chromatic-choosable only when either the cycle is odd or the cycle is even and the path is $P_1$. For $m\geq 2$, $\chi_\ell(P_m \Box C_{2n+2})=3$ and hence it is not chromatic-choosable. The graph $(K_m \lor C_{2k+1})\Box P_n$ is chromatic-choosable for $m,k,n \in \mathbb{N}$. \\

\noindent \textbf{Open problem \cite{kaul2018alon}:} \textit{For what graphs $G$, is $G \Box P_n$ chromatic-choosable?}\\

\noindent \textbf{Open problem \cite{kaul2018alon}:} \textit{When is $C_{2k+1} \Box P_n^r$ chromatic-choosable?}\\
 
\par A few other results on list-coloring of the Cartesian product of graphs are discussed in \cite{kaul2019list} and \cite{kaul2021criticality}. We do not discuss these results in detail, since they involve strong chromatic choosability and the list-color function.\\

Petrov and Gordeev \cite{petrov2022alon} have given some results for the Cartesian product involving even cycles. The first result was that $C_{2k} \Box K_n=n$, hence chromatic-choosable. The second result was that for $p,n \in \mathbb{N}$, if $(p+1)$ divides $n$ or $n \geq p(p+1)$, then $C_{2k} \Box C_n^p$ is $(p+2)$-choosable. The final result is the Cartesian product of more than two graphs. Let $0\leq m<n$ and $k_1,k_2, \dots, k_m$ be such that $\frac{1}{k_1}+\frac{1}{k_2}+\dots+\frac{1}{k_m} \leq 1$. If $G=C_{2k_1+1} \Box C_{2k_2+1} \Box \dots \Box C_{2k_m+1} \Box C_{2k_{m+1}} \Box \dots \Box C_{2k_n}$, then $G$ is $(n+1)$-choosable. \\

\par Lu et al.\cite{li2023alon} have worked on finding the Alon-Tarsi number of the Cartesian product of two cycles. And using this, they have shown that $\chi_\ell(C_{2m} \Box C_{2n+1})= AT(C_{2m}\Box C_{2n+1})= 3$ and $\chi_\ell(C_{2m}\Box C_{2n})= AT(C_{2m}\Box C_{2n})=3$. That is, the Cartesian product of an odd and an even cycle is chromatic-choosable, while the product of two even cycles is not chromatic-choosable. Even though the Alon-Tarsi number of the Cartesian product of two odd cycles is obtained as 4, the list chromatic number is still not obtained. They conjectured it as follows.

\begin{conjecture}[\cite{li2023alon}]
      $\chi_\ell(C_{2m+1}\Box C_{2n+1})=3$.
\end{conjecture}
\par The most recent work is given by Li et al.\cite{li2025alon}. They have succeeded in finding the Alon-Tarsi number of the Cartesian product of a hypercube with a tree and an even cycle. They have also pointed to the fact that in either case, the list chromatic number is strictly less than the Alon-Tarsi number. Finding the list chromatic number of these products remains open. The authors have given an open question at the end of the paper.\\

\noindent\textbf{Open problem \cite{li2025alon}:} \textit{For a given pair of graphs $G$ and $H$ with some conditions, if $\chi_\ell(G)=AT(G)$ and $\chi_\ell(H)=AT(H)$, is $\chi_\ell(G \Box H)=AT(G \Box H)$?}
    
\subsection{Lexicographic product}
\hfill \par The lexicographic product of two graphs $G=(V(G),E(G))$ and $H=(V(H),E(H))$, denoted by $G[H]$, is the graph with vertex set $V(G) \times V(H)$ and any two vertices $(u,v),(u',v')$ are adjacent if and only if $uu' \in E(G)$ or $u=u'$ and $vv' \in E(H)$. As with the other products, Hammack et al. \cite{hammack2011handbook} proved that $\chi(G[H]) \leq \chi(G)\chi(H)$. \\

\par Bal{\'a}zs and Xuding \cite{keszegh2017choosability} gave a bound for the lexicographic product of two graphs of the same order. Let $G$ and $H$ be graphs with $n$ vertices. Then $\chi_\ell(G[H]) \leq (4\Delta (G) + 2)(\chi_\ell(H) + \log n)$. They also posed a question on a slightly different bound for lexicographic product of two graphs of the same order. \\

\noindent\textbf{Open problem \cite{kim2015chromatic} :} \textit{Does there exist a constant $c$ such that, for any graph $G$ and any graph $H$ on $n$ vertices, $\chi_\ell(G[H]) \leq c \chi_\ell(G)(\chi_\ell(H) + \log n)$?}\\

\par \par Let $G$ be a graph with $V(G)=\{u_1, u_2, \dots, u_n\}$. A graph $G'$ is the \textit{inflation} of $G$ if $V(G')$ is the disjoint union of $n$ sets $V_1,V_2,\dots,V_n$ such that each $V_k$ is a clique and $x \in V_i, y \in V_j$ are adjacent in $G'$ if and only if $v_i,v_j$ is adjacent in $G$. Kostochka and Woodall \cite{kostochka2001choosability} had proved that the inflation of a graph with at most five vertices and the inflation of square of a graph with at most seven vertices are chromatic-choosable. The method employed for the proof is Hall's theorem. \\

Let $G'$ be an inflation of the graph $G$. If all the $V_k$'s are of same size, say $t$, it is called \textit{uniform inflation} and is denoted by $(G)_t$. It can be observed that this uniform inflation is the same as the lexicographic product $G[K_t]$. Prowse and Woodall \cite{prowse2003choosability}  showed that the uniform inflation of power of a cycle is chromatic-choosable.i.e., $\chi_\ell((C_n^p)_t)= \chi((C_n^p)_t)$. Since it is uniform inflation, we may restate that the lexicographic product $C_n^p[K_t]$ is chromatic-choosable. 

\section{Some regular graphs}\label{Some regular graphs}
\hfill \par Fleischner and Stiebitz \cite{fleischner1992solution} proved that for a 4-regular graph $G$ on $3n$ vertices whose edges can be decomposed into a Hamiltonian circuit and $n$ pairwise vertex-disjoint triangles, $AT(G)=3$. As $\chi(G)=3$, we have, any graph satisfying these conditions are chromatic-choosable. Later,
Gutner and Tarsi \cite{gutner2009some} gave a generalization of this result for any $(3k+1)$-regular graph on $3kn$ vertices whose edges can be decomposed into a Hamiltonian circuit and $n$ pairwise vertex-disjoint $3k$-cliques. They showed that for such a graph $G$, $\chi_\ell(G)=\chi(G)=3k$.

\section{Some variants of list-coloring}\label{Some variants of list-coloring}
\subsection{Paintability}
\hfill \par The concept of paintability was introduced in 2009 by Schauz\cite{schauz2009mr}. Let us look at a game on a given graph $G$.  Paintability is defined through a two-player game on a graph $G$, played by Lister and Painter. Let $f: V(G) \to \mathbb N$. Initially, no vertex is colored, and each vertex $v$ has $f(v)$ tokens. In each round, Lister chooses a set $S$ of uncolored vertices and removes one token from every vertex in $S$. Painter then colors an independent subset of $S$. Lister wins if some uncolored vertex has no tokens remaining; Painter wins if all vertices are colored. The graph $G$ is $f$-{\em paintable } if Painter has a winning strategy for the function $f$, and $G$ is $k$-{\em paintable } if it is $f$-paintable for the constant function $f(v)=k$. The {\em paint number } $\chi_p(G)$ is the least $k$ for which $G$ is $k$-paintable.

\par It can be seen that $\chi_\ell(G) \leq \chi_p(G)$ and the difference between can be arbitrarily large \cite{duraj2015chip}.

\subsection{DP coloring}
\hfill \par Dvo{\v{r}}{\'a}k and Postle \cite{dvovrak2018correspondence}  introduced a generalization of list-coloring called \textit{correspondence coloring}, which was later renamed as \textit{DP-coloring}. \\

\par Let $G$ be a graph. A \textit{cover} of $G$ is defined to be a pair $(L,H)$, where $L$ is a pairwise disjoint list assignment of $G$ and $H$ is a graph with $V(H)=\cup_{v \in V(G)}L(v)$ with the following conditions:
\par (a) For each $v \in V(G)$, $H[L(v)]$ is a complete graph.
\par (b) For each $uv \in E(G)$, the edge set between $L(u)$ and $L(v)$ is a matching.
\par (c) For each $uv \notin E(G)$,  the edge set between $L(u)$ and $L(v)$ is empty.\\
An \textit{$(L,H)$-coloring} of $G$ is an independent set $I$ of $V(H)$ of size $n(G)$. And $G$ is \textit{$(L,H)$-colorable} if $G$ has an $(L,H)$-coloring. The \textit{DP-chromatic number} of a graph $G$, denoted by $\chi_{DP}(G)$, is the least $k$ for which $G$ is $(L,H)$-colorable for every cover $(L,H)$, with $L$ being a $k$-list assignment.\\

\par Since DP-coloring generalizes list-coloring, we have $\chi_\ell(G) \leq \chi_{DP}(G)$.

\subsection{$(a:b)$-choosability}
\hfill \par Another early variant introduced by  Erd\H{o}s, Rubin, and Taylor \cite{erdos1979choosability}, in the same paper in which they introduced list-coloring is $(a:b)$-choosability.  

\par Let $L$ be an $a$-list assignment on $G$. If there exists a $b$-list assignment $L'$ on $G$ in such a way that 
\par (a) $L'(v) \subseteq L(v)$, for all $v \in V(G)$
\par (b) $L'(x) \cap L'(y)=\phi$ whenever $xy \in E(G)$\\
then we say that $G$ is \textit{$(a:b)$-choosable}. The \textit{$k^{th}$ choice number} of $G$, denoted by $ch_{k}(G)$, is the smallest $n$ for which $G$ is $(n:k)$-choosable. \\

\par It can be seen that $(a:1)$-choosable is nothing but $a$-choosable. Hence $\chi_\ell(G)=ch_1(G)$. Prowse and Woodall \cite{prowse2003choosability} inspected the equivalence of $(a:b)$-choosable and $a$-choosable, and conjectured as follows.

\begin{conjecture}[$(a:b)$-Choosability Equivalence Conjecture \cite{prowse2003choosability}] 
For all $a,b \in \mathbb{N}$, a graph $G$ is $(a:b)$-choosable if and only if $(G)_{b}$ is $a$-choosable.
\end{conjecture}

\begin{conjecture}[Weak $(a:b)$-Choosability Conjecture \cite{prowse2003choosability}] 
For all $a,b,k \in \mathbb{N}$, if a graph $G$ is $(a:b)$-choosable, then $G$ is $(ka:kb)$-choosable.
\end{conjecture}

\par We can find a few examples for which the conjecture holds. Tuza and Voigt \cite{tuza1996every} proved that every 2-choosable graph is $(2m:m)$-choosable for every $m \in \mathbb{N}$. Gutner and Tarsi \cite{gutner2009some} proved that every planar bipartite graph is $(3m:m)$-choosable for every $m \in \mathbb{N}$ and every $d$-degenerate graph is $((d+1)m:m)$-choosable for every $m \in \mathbb{N}$. The conjecture was disproved by Dvo{\v{r}}{\'a}k et al.\cite{dvovrak20184} by producing a 4-choosable graph which is not (8:2)-choosable.

\subsection{Fractional coloring}
\hfill \par The fractional coloring of a graph $G$ was introduced by Hilton et al.\cite{hilton19735} in the name multicoloring. It is a variant of $(a:b)$-choosability.\\

\par Let $L=\{1,2,\dots,a\}$ be the set of colors available for $G$. If there exists a $b$-list assignment $L'$ on $G$ in such a way that 
\par (a) $L'(v) \subseteq L(v)$, for all $v \in V(G)$
\par (b) $L'(x) \cap L'(y)=\phi$ whenever $xy \in E(G)$\\
then we say that $G$ has a \textit{fractional $(a:b)$-coloring}. The \textit{fractional chromatic number} of $G$, denoted by $\chi_f(G)$, is the infimum of the fractions $\frac{a}{b}$.\\

\par It can be observed that the fractional $(a,b)$-coloring on a graph is merely $(a:b)$-choosability with the additional condition that the $a$-list assignment is identical for all vertices.

\section{Conclusion}\label{Conclusion}
\hfill \par Through the survey, we have observed how list-coloring has evolved over the years and how it has developed for different classes of graphs and different operations on graphs. \\

\par First we considered the embedded graphs, we restricted our analysis to planar graphs and toroidal graphs. All planar graphs are 5-choosable and all toroidal graphs are 7-choosable. Any further tightening of the bound for list-coloring is influenced by forbidding few structures like cycles and intersecting or adjacent cycles. When there are conditions under which both planar and toroidal graphs behave alike while list-coloring, there are even more factors like girth and distance between triangles that play a crucial role for planar graphs. We have seen that the list-coloring of planar graphs have witnessed some good constructions of counterexamples - two different 4-colorable non-4-choosable planar graphs, a 3-colorable non-4-choosable planar graph, non-3-choosable planar graphs without $C_4,C_5$ and intersecting triangles and two different non-3-choosable planar graphs with girth 4.\\

\par List-coloring of complete bipartite and  multipartite graphs have been studied in the literature. The list chromatic number is found for $K_{p,q}$ for $p=1,2,3,4$ and $q \geq p$. Whereas $\chi_\ell(K_{p,q})$ is not known for a certain range of $q$ for $p=5,6,7$, it is unknown for  $p>7$ and $p \leq q \leq (p-1)^{p-1}-(p-2)^{p-1}$ except that they are $p$-choosable.  Even though the chromatic-choosable complete multipartite graphs have not been characterized, the Ohba conjecture, being proved true, implies that if the total number of vertices is at most one more than its list-coloring. Also the list chromatic number of a few other complete multipartite graphs are shown.\\

\par Claw-free graphs, another important class, have also gained importance in the problem of list-coloring. There have been few results on the list-coloring of claw-free graphs with small $\Delta$ and claw-free perfect graphs with small clique number. As the coloring of line graph of a graph can be equalized by edge coloring the graph itself, we had to deviate into list edge coloring and the results on it.\\

\par The list-coloring of powers of graphs has emerged as an important problem in chromatic-choosability. Considerable progress has been made for squares of planar graphs, especially for subcubic planar graphs. One recurring theme is that, under suitable restrictions on girth and maximum degree, the list chromatic number of the square of a planar graph becomes more tractable and, in many cases, admits sharp or near-sharp bounds. Although maximum degree is fixed for subcubic planar graph, girth still plays a role. Among the basic graph classes, paths and cycles remain the principal examples whose powers are known to be chromatic-choosable for every power.\\

\par Finding the coloring parameter for different graph operations in terms of the parameters of the component graphs is one classic challenge. Ohba introduced the condition under which the join of two graphs becomes chromatic-choosable. Analyzing the list-coloring problem for the product, the significant progress is limited to the Cartesian product and the Lexicographic product. When the bound for the list chromatic number of the Cartesian product of $k$ graphs is found, for the Lexicographic product, the bound is only for the product of two graphs of same order.\\

The results surveyed here show that chromatic choosability is governed by a diverse range of structural and methodological phenomena. Degeneracy and Hall-type arguments provide some of the earliest tools, while the Alon–Tarsi method, kernel-perfect orientations, boundary-list induction, and discharging play central roles in more refined results. Several major conjectures, such as the List Square Coloring Conjecture and its planar version, have been disproved by explicit counterexamples, while others remain open or are known only under additional structural assumptions. The remaining problems suggest several promising directions, including claw-free graphs, line graphs, graph powers, complete multipartite graphs, and the behavior of list coloring under graph operations.

\printbibliography[heading=bibintoc, title={Reference}]
\end{document}